\documentclass{amsart}	
\usepackage{amsmath}
\usepackage{amssymb}
\usepackage{amsfonts}
\usepackage{amsthm}

\usepackage{tcolorbox}
\tcbuselibrary{breakable, skins, theorems}
\tcbuselibrary{breakable}

\usepackage{mathtools, enumerate, enumitem, refcount, xparse}
\mathtoolsset{showonlyrefs = true}

\usepackage{hyperref}
\hypersetup{
	setpagesize=false,
	bookmarksnumbered=true,
	bookmarksopen=true,
	colorlinks=true,
	linkcolor=blue,
	citecolor=blue,
}

\makeatletter
\@addtoreset{equation}{section}

\makeatother

\newenvironment{eqenumerate}
{
	\begin{enumerate}[ref=\thesection.\theenumi, leftmargin=!, align=left]
		\setcounter{enumi}{\value{equation}}
		
	}
	{
		\setcounter{equation}{\value{enumi}}
	\end{enumerate}
}

\NewDocumentCommand\eqitem{ o }{
	\setcounter{enumi}{\value{equation}}
	\IfValueTF{#1}
	{\item[#1]}
	{\item}
	\setcounter{equation}{\value{enumi}}
}

\theoremstyle{plain}
\newtheorem{thm}{Theorem}[section]
\newtheorem{prop}[thm]{Proposition}
\newtheorem{lem}[thm]{Lemma}

\theoremstyle{definition}
\newtheorem{dfn}[thm]{Definition}

\theoremstyle{remark}
\newtheorem*{rem}{Remark}

\newcommand{\R}{\mathbb{R}}

\newcommand{\ria}{\rightarrow}

\newcommand{\lan}{\langle}
\newcommand{\ran}{\rangle}

\newcommand{\lam}{\lambda}
\newcommand{\Ag}{\mathcal{A}_{\gamma}}
\newcommand{\Az}{\mathcal{A}_{0}}
\newcommand{\hsl}{H^{s/2}\times L^2}
\newcommand{\hsld}{H^{s/2}(\R^d)\times L^2 (\R^d)}
\newcommand{\hsh}{H^{s}\times H^{s/2}}
\newcommand{\hshd}{H^s (\R^d)\times H^{s/2}(\R^d)}
\newcommand{\hsho}{H^s (\R) \times H^{s/2}(\R)}
\newcommand{\hslo}{H^{s/2} (\R) \times L^2 (\R)}
\newcommand{\F}{\mathcal{F}}

\DeclareMathOperator{\supp}{supp}
\DeclareMathOperator*{\essinf}{ess\,inf}
\renewcommand{\Re}{\operatorname{Re}}

\usepackage[abbrev]{amsrefs}
\BibSpec{arXiv}{%
	+{}{\PrintAuthors}{author}
	+{,}{ \textit}{title}
	+{,}{ arXiv:}{eprint}
	+{}{ \parenthesize}{date}
	+{.}{}{transition}
}

\title[Energy decay of damped Klein--Gordon equations]{Equivalence between the energy decay of fractional damped Klein--Gordon equations and geometric conditions for damping coefficients}

\author[K. Inami]{Kotaro Inami}
\address[Kotaro Inami]{Graduate School of Mathematics, Nagoya University, Furocho, Chikusa-ku, Nagoya, Aichi, 464-8602, Japan}
\email{m21010t@math.nagoya-u.ac.jp}
\thanks{The first author was supported by JST SPRING Grant Number JPMJSP2125
, the Interdisciplinary Frontier Next-Generation Researcher Program of the Tokai Higher Education and Research System
.}

\author[S. Suzuki]{Soichiro Suzuki}
\address[Soichiro Suzuki]{Department of Mathematics, Chuo University, 1-13-27, Kasuga, Bunkyo-ku, Tokyo, 112-8551, Japan}
\email{soichiro.suzuki.m18020a@gmail.com}
\thanks{The second author was supported by Japan Society for the Promotion of Science (JSPS) KAKENHI Grant Number JP20J21771 and JP23KJ1939.}


\begin{document}
	\begin{abstract}
We consider damped $s$-fractional Klein--Gordon equations on $\mathbb{R}^d$, where $s$ denotes the order of the fractional Laplacian. 
In the one-dimensional case $d = 1$, Green (2020) established that the exponential decay for $s \geq 2$ and the polynomial decay of order $s/(4-2s)$ hold if and only if the damping coefficient function satisfies the so-called geometric control condition. In this note, we show that the $o(1)$ energy decay is also equivalent to these conditions in the case $d=1$. 
Furthermore, we extend this result to the higher-dimensional case: the logarithmic decay, the $o(1)$ decay, and the thickness of the damping coefficient are equivalent for $s \geq 2$. 
In addition, we also prove that the exponential decay holds for $0 < s < 2$ if and only if the damping coefficient function has a positive lower bound, so in particular, we cannot expect the exponential decay under the geometric control condition.
\end{abstract}

	\maketitle
	
	\section{Introduction} \label{section:intro}
	We consider the following fractional damped Klein--Gordon equations on $\R^d$:
\begin{equation}
u_{tt}(t,x) + \gamma(x)u_{t}(t,x) + (-\Delta + 1)^{s/2}u(t,x) = 0, \quad (t,x)\in \R_{\geq 0}\times \R^d, 
\label{eq:DWE}
\end{equation}
where $s > 0$, and $0 \leq \gamma \in L^\infty(\R^d)$. Here we note that $\gamma u_{t}$ represents the damping force and the operator $(-\Delta + 1)^{s/2}$ is defined by the Fourier transform on $L^2 (\R^d)$; 
\begin{equation*}
(-\Delta + 1)^{s/2}u \coloneqq \mathcal{F}^{-1}(|\xi|^2 + 1)^{s/2}\mathcal{F}u,\quad \xi\in\R^{d}. 
\end{equation*}
We recast the equation \eqref{eq:DWE} as an abstract first-order equation for $U = (u,u_{t})$:
\begin{align}
U_{t} = \Ag
U, \quad \Ag = \begin{pmatrix}
0 & I \\
-(-\Delta + 1)^{s/2} & -\gamma(x)
\end{pmatrix}, \label{eq:ode}
\end{align}
then $\Ag$ generates a $C_{0}$-semigroup $(e^{t A_{\gamma}})_{t \geq 0}$ on $\hsld$ (see \cite{green2020energy}). Here the Sobolev space $H^{r}(\R^d)$ is defined by
\begin{align*}
H^r (\R^d) \coloneqq \left\{u\in L^2 (\R^d) : \|u\|^{2}_{H^r} = \int_{\R^d}(|\xi|^2 + 1)^r |\F u(\xi)|^2 \,d\xi < \infty\right\}. 
\end{align*}
In this paper, we discuss the decay rate of the energy
\begin{equation*}
E(t) \coloneqq \|e^{t\Ag}(u(0),u_{t}(0))\|_{H^{s/2}\times L^2} = \left(\int_{\R^d} ( |(-\Delta + 1)^{s/4}u(t,x)|^2 + |u_{t}(t,x)|^2 ) \,dx\right)^{1/2}
\label{energy}. 
\end{equation*}
By standard calculus, we have $E(t) = E(0)$ if $\gamma \equiv 0$ and the exponential energy decay if $\gamma \equiv C > 0$. In recent works, the intermediate case, that is, the case that $\gamma = 0$ on a large set is studied: 
\begin{dfn}
We say that $\Omega\subset\R^d$ satisfies the \textit{Geometric Control Condition} $\mathrm{(GCC)}$ if there exist $L>0$ and $0 < c \leq 1$ such that for any line segments $l\in\R^d$ of length $L$, the inequality 
\begin{align*}
\mathcal{H}^{1}(\Omega\cap l) \geq cL
\end{align*}
holds,
where $\mathcal{H}^{1}$ denotes the one-dimensional Hausdorff measure. 
\end{dfn}
Burq and Joly \cite{burq2016exponential} proved that if $\gamma$ is uniformly continuous and $\{\gamma\geq\varepsilon\}$ satisfies (GCC) for some $\varepsilon > 0$, then we have the exponential energy decay in the non-fractional case $s = 2$. 
After that, Malhi and Stanislavova \cite{malhi2018energy} pointed out that (GCC) is also necessary for the exponential decay in the one-dimensional case $d = 1$:
\begin{thm}[{\cite[Theorem 1]{malhi2018energy}}] \label{Malhi1}
Let $d = 1$, let $s = 2$, and let $0 \leq \gamma \in L^{\infty}(\R)$ be continuous. Then the following conditions are equivalent: 
\begin{eqenumerate}
\eqitem\label{item:gcc} There exists $\varepsilon > 0$ such that the upper level set $\{\gamma\geq\varepsilon\}$ satisfies $\mathrm{(GCC)}$. 
\eqitem\label{item:expenergy} There exist $C,\omega > 0$ such that whenever $(u(0),u_{t}(0))\in H^1(\R)\times L^2(\R)$,
\begin{align*}
E(t) \leq C\exp(-\omega t) E(0) 
\end{align*}
holds for any $t \geq 0$.%
\footnote{To be precise, the exponential decay estimate given in \cite[Theorem 1]{malhi2018energy} is a little weaker: $E(t) \leq C\exp(-\omega t) \|(u(0),u_{t}(0))\|_{H^2\times H^1}$. However, this is because the Gearhart--Pr\"uss theorem in their paper (\cite[Theorem 2]{malhi2018energy}) is stated incorrectly. Using the theorem correctly (see Theorem \ref{thm:gearhart-Pruss}), one can obtain the exponential decay estimate $E(t) \leq C\exp(-\omega t) E(0)$ as in \eqref{item:expenergy}.}
\eqitem\label{item:energydecay s=2}
$\underset{t\ria +\infty}{\lim}\|e^{t A_{\gamma}}\|_{H^2\times H^{1}\ria H^{1}\times L^2} = 0$.
\end{eqenumerate}
\end{thm}
Note that for $0 \leq \gamma \in L^{\infty}(\R)$, the condition \eqref{item:gcc} is also equivalent to that there exists $R > 0$ such that 
\begin{equation*}
\underset{a\in \R}{\inf}\int^{a + R}_{a - R}\gamma(x) \,dx > 0. 
\end{equation*}
In another paper \cite{malhi2020energy}, Malhi and Stanislavova introduced the fractional equation \eqref{eq:DWE} and showed that if $\gamma$ is periodic, continuous and not identically zero, then we have the exponential decay for any $s \geq 2$ and the polynomial decay of order $s/(4-2s)$ for any $0 < s < 2$ in the case $d = 1$.
\begin{rem}
Nonzero periodic functions satisfy (GCC) in the case $d = 1$, but it is not true in the higher-dimensional case $d \geq 2$. 
Wunsch \cite{wunsch2017periodic} showed that continuous periodic damping gives the polynomial energy decay of order $1/2$ for the non-fractional equation in the case $d \geq 2$.
In addition, recently another proof and an extension to fractional equations of Wunsch's result were obtained by T\"aufer \cite{Taufer2022} and Suzuki \cite{Suzuki2022}, respectively.
Note that these results for periodic damping are established by reducing to estimates on the torus $\mathbb{T}^d$. 
Indeed, there are numerous studies on bounded domains; see references in \cite{burq2016exponential} and \cite{CST2020}, for example.
\end{rem}
Green \cite{green2020energy} improved results of Malhi and Stanislavova as follows: 
\begin{thm}[{\cite[Theorem 1]{green2020energy}}]\label{Green1}
Let $d = 1$, let $s > 0$ and let $0 \leq \gamma \in L^{\infty}(\R)$. Then the following conditions are equivalent: 
\begin{eqenumerate}
\eqitem[\eqref{item:gcc}] There exists $\varepsilon > 0$ such that the upper level set $\{\gamma\geq\varepsilon\}$ satisfies $\mathrm{(GCC)}$. 
\eqitem\label{item:exppolenergy} There exist $C,\omega > 0$ such that whenever $(u(0),u_{t}(0))\in\hsho$, 
\begin{align*}
E(t) \leq 
\begin{cases}
(1 + t)^{-\frac{s}{4 - 2s}}\|(u(0),u_{t}(0))\|_{\hsh}\quad &\text{if} \quad 0 < s < 2,\\
C\exp(-\omega t )E(0)\quad &\text{if} \quad s \geq 2
\end{cases}
\end{align*}
holds for any $t \geq 0$. 
\end{eqenumerate}
\end{thm}
In comparison with the result of \cite{malhi2020energy}, which states that \eqref{item:exppolenergy} holds if $\gamma$ is periodic, continuous and not identically zero,
Theorem \ref{Green1} refines this result by giving a necessary and sufficient condition for \eqref{item:exppolenergy}. 
Furthermore, Theorem \ref{Green1} also improves the $\eqref{item:gcc} \iff \eqref{item:expenergy}$ part of Theorem \ref{Malhi1} by extending it to fractional equations and removing the continuity of $\gamma$, 
but on the other hand, it lacks the $\eqref{item:energydecay s=2} \implies \eqref{item:gcc}, \eqref{item:expenergy}$ part. 
One of our goal is to recover this part for fractional equations: 
\begin{thm}\label{thm1}
Let $d = 1$, let $s > 0$, and let $0\leq \gamma \in L^{\infty}(\R)$. Then the following conditions are equivalent: 
\begin{eqenumerate}
\eqitem[\eqref{item:gcc}] There exists $\varepsilon > 0$ such that the upper level set $\{\gamma\geq\varepsilon\}$ satisfies $\mathrm{(GCC)}$. 
\eqitem[\eqref{item:exppolenergy}] There exist $C,\omega > 0$ such that whenever $(u(0),u_{t}(0))\in\hsho$,
\begin{align*}
E(t) \leq \begin{cases}
C(1 + t)^{\frac{-s}{4 - 2s}}\|(u(0),u_{t}(0))\|_{\hsh} \quad &\text{if} \quad 0 < s < 2, \\
C\exp(-\omega t)E(0) \quad &\text{if} \quad s\geq 2
\end{cases}
\end{align*}
holds for any $t \geq 0$. 
\eqitem \label{item:energydecay} $\underset{t\ria +\infty}{\lim}\|e^{t A_{\gamma}}\|_{H^s\times H^{s/2}\ria H^{s/2}\times L^2} = 0$.
\end{eqenumerate}
\end{thm}
We also give the following result, which says that we cannot expect the exponential decay for $0 < s < 2$ under (GCC). 
\begin{thm}\label{thm:exponentialdecay}
Let $ d\geq 1$, let $0 < s < 2$, and let $0\leq\gamma\in L^\infty (\R^d)$. Then there exist $C,\omega > 0$ such that whenever $(u(0), u_{t}(0))\in\hsld$, 
\begin{align}
E(t) \leq C\exp(-\omega t)E(0)
\end{align}
holds for any $t \geq 0$
if and only if $\essinf_{\R^d}\gamma>0$. 
\end{thm}
Note that the ``if'' part easily follows by reducing to the constant damping case, so we will prove the ``only if'' part.
Furthermore, we extend Theorem \ref{thm1} to the higher-dimensional case $d\geq 2$ using a notion of \textit{thickness}, which is equivalent to (GCC) in the case $d = 1$:
\begin{dfn}
We say that a set $\Omega \subset\R^d$ is \textit{thick} if there exists $R > 0$ such that
\begin{align}
\underset{a\in\R^d}{\inf}m_{d}(\Omega\cap (a + [-R, R]^d)) > 0 \label{dfnthick}
\end{align}
holds, where $m_{d}$ denotes the $d$-dimensional Lebesgue measure. 
\end{dfn}
Then we have the following result:
\begin{thm}\label{thm:thick}
Let $d \geq 2$, let $s \geq 2$, and let $0\leq \gamma\in L^\infty (\R^d)$. Then the following conditions are equivalent: 
\begin{eqenumerate}
\eqitem\label{item:thickg}
There exists $\varepsilon > 0$ such that the upper level set $\{\gamma \geq \varepsilon\}$ is thick. 
\eqitem\label{item:logenergy}
There exists $C > 0$ such that whenever $(u(0),u_{t}(0))\in\hshd$,
\begin{equation*}
E(t) \leq \frac{C}{\log(e + t)}\|(u(0),u_{t}(0))\|_{\hsh}
\end{equation*}
holds for any $t \geq 0$. 
\eqitem\label{item:denergydecay}
$\underset{t\ria +\infty}{\lim}\|e^{t A_{\gamma}}\|_{H^s\times H^{s/2}\ria H^{s/2}\times L^2} = 0$.
\end{eqenumerate}
\end{thm}
The implication $\eqref{item:thickg} \implies \eqref{item:logenergy}$ is a generalization of the result given by Burq and Joly \cite{burq2016exponential}. They established \eqref{item:logenergy} under the so-called network control condition, which is stronger than \eqref{item:thickg}.
Also, similarly to the case $d = 1$, the condition \eqref{item:thickg} is equivalent to that there exists $R > 0$ such that 
\begin{equation*}
\underset{a\in \R}{\inf}\int_{ a +[-R, R]^d } \gamma(x) \,dx > 0. 
\end{equation*}

Finally, we explain the organization of this paper. In Sections \ref{section:1}, \ref{section:2}, and \ref{section:3}, we will give proofs of Theorems \ref{thm1}, \ref{thm:exponentialdecay}, and \ref{thm:thick}, respectively. To prove these theorems, we use a kind of uncertainty principle and results of the $C_{0}$ semigroup theory.

	\section{Proof of Theorem \ref{thm1}} \label{section:1}
	To prove this theorem, we use the following result by Batty, Borichev, and Tomilov \cite{batty2016lp}:
\begin{thm}[{\cite[Theorem 1.4]{batty2016lp}}]
Let $A$ be a generator of a bounded $C_{0}$-semigroup $(e^{tA})_{t\geq 0}$ on a Banach space $X$, and $\lambda \in\rho (A)$. Then the following are equivalent: 
\begin{eqenumerate}
\eqitem $\sigma(A) \cap i\R = \emptyset$, 
\eqitem $\lim_{t\ria\infty}\|e^{tA}(\lambda - A)^{-1}\|_{\mathcal{B}(X)} = 0$.
\end{eqenumerate}
\end{thm}
In the case $A = \Ag$, for $\lambda \in \rho(\Ag)$, the map $(\lambda - \Ag)^{-1} : \hslo \ria \hsho$ is surjective. Thus, we have:
\begin{lem}[{\cite[Corollary 2]{malhi2018energy}}]\label{malhiwave}
For the semigroup $e^{t\Ag}$ of the Cauchy problem $\eqref{eq:ode}$, the following are equivalent: 
\begin{eqenumerate}
\eqitem $\sigma(\Ag) \cap i\R = \emptyset$, 
\eqitem $\lim_{t\ria\infty} \|e^{t\Ag}\|_{\hsh\ria\hsl} = 0$. 
\end{eqenumerate}
\end{lem}
\begin{proof}[Proof of Theorem \ref{thm1}. ]
It is enough to show that $\eqref{item:energydecay} \implies \eqref{item:gcc}$, since $\eqref{item:gcc} \iff \eqref{item:exppolenergy}$ is already known by Green \cite{green2020energy} (Theorem \ref{Green1}) and $\eqref{item:exppolenergy} \implies \eqref{item:energydecay}$ is trivial.
Suppose that \eqref{item:energydecay} holds, that is, $\lim_{t \to +\infty} \|e^{t A_{\gamma}}\|_{H^s\times H^{s/2}\ria H^{s/2}\times L^2} = 0$.
By Lemma \ref{malhiwave}, we have $i\R\subset \rho(\mathcal{A}_{\gamma})$. 
This implies that for each $\lambda \in \R$, there exists some $c_0 > 0$ such that
\begin{equation}
c_0 \|U\|_{H^{s/2}\times L^2} \leq \|(\Ag - i\lam I)U\|_{H^{s/2}\times L^2}
\end{equation}
holds for any $U\in H^{s}(\R)\times H^{s/2}(\R)$. Letting $u\in L^2 (\R^d)$ and $U = ((-\Delta + 1)^{-s/4}u, iu)$, we obtain
\begin{align*}
2c_0 \|u\|^{2}_{L^2} &\leq \|((-\partial_{xx} + 1)^{s/4} - \lam)u\|^{2}_{L^2} + \|( (-\partial_{xx} + 1)^{s/4}  - \lambda + i\gamma)u\|^{2}_{L^2}\\
&\leq 3\|((-\partial_{xx} + 1)^{s/4} - \lam)u\|^{2}_{L^2} + 2\|\gamma u\|^{2}_{L^2}.
\end{align*}
Now we consider the case $\lam = 1$. Let $u\in H^{s/2}(\R)$ satisfy $\supp\widehat{u} \subset [-D,D]$ for some $D>0$, which is chosen later. For such $u$, we have
\begin{align*}
\|((-\partial_{xx} + 1)^{s/4} - 1)u\|^{2}_{L^2} &= \int^{D}_{-D}\left[(|\xi|^2 + 1)^{s/4} -1\right]^2 |\widehat{u}(\xi)|^2 d\xi\\
&\leq \left[(D^2 + 1)^{s/4} -1\right]^{2}\|u\|^{2}_{L^2}.
\end{align*}
Hence, taking $D>0$ small enough, we get some $c > 0$ such that
\begin{equation}
c \|u\|_{L^2} \leq \|\gamma u\|_{L^2}\label{eq:neqgamma}
\end{equation}
holds for any $u\in H^{s/2}(\R)$ satisfying $\supp\widehat{u}\subset [-D,D]$.
Fix $f \in \mathcal{S}(\R) \setminus \{ 0 \}$ such that $\supp \widehat{f} \subset [-D, D]$
and write $f_a(x) \coloneqq f(x-a)$ for each $a\in\R$, so that $\widehat{f_a}(\xi) = e^{ia\xi} \widehat{f}(\xi)$. 
Then, for each $a \in \R$ and $R > 0$, we have
\begin{equation*}
0 < c \|f\|_{L^2} = c \| f_a \|_{L^2} \leq \| \gamma f_a \|_{L^2} = \left(\int_{[a - R,a + R]} + \int_{[a-R,a+R]^{c}} \right)|\gamma(x) {f}_{a}(x)|^2 \,dx.
\end{equation*}
The second integral goes to $0$ as $R \to + \infty$ since $\gamma$ is bounded and $|{f}_{a}|^2$ is integrable, and this convergence is uniform with respect to $a$. 
Furthermore, for the first integral, we have
\begin{equation*}
\int^{a+R}_{a-R}|\gamma(x) {f}_{a}(x)|^{2} \,dx \leq \|\gamma\|_{L^{\infty}}\|f\|^{2}_{L^{\infty}}\int^{a+R}_{a-R}\gamma(x) \,dx ,
\end{equation*}
since $\gamma$ and $f$ are bounded and $\|f_a\|_{L^\infty} = \|f\|_{L^\infty}$. Thus, there exists $R>0$ such that
\begin{equation*}
\inf_{a \in \R} \int^{a + R}_{a - R}\gamma(x) \,dx > 0
\end{equation*}
holds, 
which is equivalent to \eqref{item:gcc}.
\end{proof}

	\section{Proof of Theorem \ref{thm:exponentialdecay}} \label{section:2}
	This section is based on the proof of Theorem 2 in Green \cite{green2020energy}. To prove this theorem, we use the classical semigroup result by Gearhart, Pr\"{u}ss, and Huang:
\begin{thm}[Gearhart--Pr\"{u}ss--Huang] \label{thm:gearhart-Pruss}
Let $X$ be a complex Hilbert space and let $(e^{tA})_{t\geq 0}$ be a bounded $C_{0}$-semigroup on $X$ with infinitesimal generator $A$. Then there exist $C,\omega > 0$ such that 
\begin{align*}
\|e^{tA}\| \leq C\exp(-\omega t)
\end{align*}
holds for any $t \geq 0$
if and only if $i\R\subset \rho(A)$ and $\sup_{\lam\in\R}\|(i\lam - A)^{-1}\|_{\mathcal{B}(X)} < \infty$. 
\end{thm}

\begin{proof}[Proof of Theorem \ref{thm:exponentialdecay}]
We will prove the contraposition of the ``only if'' part of Theorem \ref{thm:exponentialdecay}, that is, 
if the energy decays exponentially and $\essinf_{x\in\R^d} \gamma(x) = 0$ holds, then $s \geq 2$. 
By the Gearhart--Pr\"{u}ss--Huang theorem and the exponential decay, there exists $c_0 > 0$ such that
\begin{align*}
c_0 \|U\|^2_{\hsl} \leq \|(\Ag - i \lam I)U\|^{2}_{\hsl}
\end{align*}
holds for any $U\in\hsld$ and any $\lam\in\R$. 
Letting $u\in L^2 (\R^d)$ and $U = ((-\Delta + 1)^{-s/4}u, iu)$, we obtain
\begin{align*}
2c_0 \|u\|^{2}_{L^2} &\leq \|( (-\Delta + 1)^{s/4} -\lam )u\|^2_{L^2} + \|((-\Delta + 1)^{s/4}  - \lambda + i\gamma)u\|^2_{L^2}\\
&\leq 3\|(-\Delta + 1)^{s/4} - \lam\|^{2}_{L^2} + 2\|\gamma u\|^{2}_{L^2}.
\end{align*}
Now let $u\in L^2 (\R^d)$ satisfy 
\begin{equation}
\supp \widehat{u} \subset \{\xi\in\R^d : |(|\xi|^2 + 1)^{s/4} - \lam| \leq K\} \eqqcolon A_{\lam}(K)
\end{equation}
for some $K$, which is chosen later.
For such $u$, we have
\begin{align*}
\|((-\Delta + 1)^{s/4} -\lam)u\|^2_{L^2} &= \int_{A_{\lam}(K)} [(|\xi|^2 + 1)^{s/4} - \lam]^2 |\widehat{u}(\xi)|^2 \,d\xi\\
&\leq K^2 \|u\|^{2}_{L^2}. 
\end{align*}
Hence, taking $K > 0$ small enough, we get some $c>0$ such that
\begin{equation}
c \|u\|^{2}_{L^2} \leq \|\gamma u\|^{2}_{L^2} \label{gamma}
\end{equation}
holds for any $u\in L^2 (\R^d)$ satisfying $\supp \widehat{u} \subset A_{\lam}(K)$ with some $\lambda \in \R$.

We prove $s \geq 2$ by contradiction. 
Assume that $s < 2$. In this case, the thickness of the annulus $A_{\lam}(K)$ is unbounded with respect to $\lam$:
\begin{equation*}
\underset{\lam\ria\infty}{\lim}\left|\sqrt{(\lam + K)^{4/s} - 1} - \sqrt{(\lam - K)^{4/s} - 1}\right| = \underset{\lam\ria\infty}{\lim}\frac{\lam^{4/s - 1}}{\lam^{2/s}} = \infty. 
\end{equation*}
Thus, the inequality \eqref{gamma} holds for any $u\in L^2 (\R^d)$ such that $\supp \widehat{u}$ is compact. 
To see this, notice that there exist $a \in \R^d$ and $\lambda \in \R$ satisfying $a + \supp \widehat{u} \subset A_\lambda(K)$ for such $u$. 
Therefore, letting $u_a(x) \coloneqq e^{i a \cdot x} u(x)$, we have
\begin{equation*}
c \|u\|^{2}_{L^2} = c \|u_a\|^{2}_{L^2} \leq \|\gamma u_a\|^{2}_{L^2} = \|\gamma u\|^{2}_{L^2}
\end{equation*}
since $\supp \widehat{ u_a } = a + \supp \widehat{u} \subset A_\lambda(K)$.

Now note that $E_{\varepsilon} \coloneqq \{x\in\R^d : \gamma(x) < \varepsilon\}$ has a positive measure for any $\varepsilon > 0$, since $\essinf_{x \in \R^d} \gamma(x) = 0$.
For each $\varepsilon > 0$, we take a subset $F_{\varepsilon} \subset E_{\varepsilon}$ such that $0 < m_{d}(F_{\varepsilon}) < \infty$. Take $R,\varepsilon > 0$ arbitrarily and set
\begin{align*}
f_{\varepsilon} \coloneqq \chi_{F_{\varepsilon}}/\sqrt{m_{d}(f_{\varepsilon})}, \quad g_{R,\varepsilon} \coloneqq \F^{-1}\chi_{B(0,R)}\F f_{\varepsilon} ,
\end{align*}
where $\chi_{\Omega}$ denotes the indicator function of $\Omega\subset\R^d$. 
By the definition, we have $\supp \widehat{g_{R,\varepsilon}} \subset B(0, R)$ and $g_{R,\varepsilon}\ria f_{\varepsilon}$ as $R\ria \infty$ in $L^2(\R^d)$. 
Therefore, applying the inequality \eqref{gamma} to $g_{R,\varepsilon}$, we get\begin{align*}
c \|g_{R,\varepsilon}\|_{L^2} &\leq \|\gamma g_{R,\varepsilon}\|_{L^2}\\
&\leq \|\gamma f_{\varepsilon}\|_{L^2} + \|\gamma(g_{R,\varepsilon} - f_{\varepsilon})\|_{L^2}\\
&= \left(\frac{1}{m_{d}(F_{\varepsilon})}\int_{F_{\varepsilon}}|\gamma(x)|^2 dx\right)^{1/2} + \|\gamma(g_{R,\varepsilon} - f_{\varepsilon})\|_{L^2}\\
&\leq \varepsilon + \|\gamma(g_{R,\varepsilon} - f_{\varepsilon})\|_{L^2}. 
\end{align*}Taking the limit as $R \to + \infty$, we obtain
\begin{align*}
0 < c = c\|f_{\varepsilon}\|_{L^2} \leq \varepsilon.
\end{align*}
This is a contradiction since $\varepsilon > 0$ is arbitrary. 
\end{proof}

	\section{Proof of Theorem \ref{thm:thick}} \label{section:3}
	The proof of $\eqref{item:denergydecay} \implies \eqref{item:thickg}$ is similar to that of $\eqref{item:energydecay} \implies \eqref{item:gcc}$ in Section \ref{section:1}, and the implication $\eqref{item:logenergy} \implies \eqref{item:denergydecay}$ is trivial. 
Therefore, we will show that $\eqref{item:thickg} \implies \eqref{item:logenergy}$. 
We use a kind of the uncertainty principle to obtain a certain resolvent estimate for the fractional Laplacian:
\begin{thm}[{\cite[Theorem 3]{kovrijkine2001some}}]\label{thm:kov}
Let $\Omega\subset\R^d$ be thick. Then there exist a constant $C > 0$ such that for each $R > 0$, the inequality
\begin{align}
\|f\|_{L^2 (\R^d)} \leq C \exp(CR)\|f\|_{L^{2}(\Omega)}
\end{align}
holds for any $f \in L^2 (\R^d)$ satisfying $\supp \widehat{f} \subset B(0,R)$.
\end{thm}
In order to obtain the logarithmic energy decay \eqref{item:logenergy}, we use the following result. 
\begin{thm}[{\cite[Theorem 5.1]{burq2016exponential}}] \label{thm:burq}
Let $A$ be a maximal dissipative operator (and hence generate the $C^{0}$-semigroup of contractions $(e^{tA})_{t \geq 0}$) in a Hilbert space $X$. Assume that $i\R \subset \rho(A)$ and there exists $C > 0$ such that
\begin{align}
\ \|(A - i\lam I)^{-1}\|_{\mathcal{B}(X)} \leq Ce^{C|\lam|}
\end{align}
holds for any $\lambda \in \R$.
Then, for each $k > 0$, there exists $C_{k} > 0$ such that
\begin{align}
\|e^{tA}(I - A)^{-k}\|_{\mathcal{B}(X)} \leq \frac{C_{k}}{( \log(e + t) )^k}
\end{align}
holds for any $t \geq 0$.
\end{thm}
\subsection{Resolvent estimate}
The proof of these propositions are based on \cite{green2020energy}. 
\begin{prop} \label{prop:resovent}
Let $s \geq 1$ and $\Omega\subset\R^d$ be thick. Then there exist $C,c>0$ such that for all $f\in L^2 (\R^d)$ and all $\lam \geq 0$, 
\begin{align}
c\exp(-C\lam)\|f\|^{2}_{L^{2}(\R^d)} \leq \|((-\Delta + 1)^{s/2} - \lam)f\|^{2}_{L^{2}(\R^d)} + \|f\|^{2}_{L^{2}(\Omega)}. 
\end{align}
\end{prop}
\begin{proof}[Proof of Proposition \ref{prop:resovent}]
Let $A_{\lam} \coloneqq \{\xi\in\R^d : |(|\xi|^2 + 1)^{1/2} - \lam^{1/s}|\leq 1\}$. Since $A_{\lam}\subset B(0,\lam + 2)$ and $\Omega$ is thick, Theorem \ref{thm:kov} implies that there exists $C > 0$ such that 
\begin{equation} \label{eq:str anni}
\|f\|_{L^2 (\R^d)} \leq C\exp(C\lam)\|f\|_{L^{2}(\Omega)}
\end{equation}
holds for any $\lambda \geq 0$ and any $f\in L^2 (\R^d)$ satisfying $\supp \widehat{f}\subset A_{\lam}$.
Next, we set a projection $P_{\lam} \coloneqq \F^{-1}\chi_{A_{\lam}}\F$, where $\chi_{A_{\lam}}$ denotes the indicator function of $A_{\lam}$. 
Then, since $P_{\lam}f$ satisfies the inequality \eqref{eq:str anni} for each $f\in L^2 (\R^d)$, we obtain
\begin{align*}
\|f\|^{2}_{L^{2}(\R^d)} &= \|P_{\lam}f\|^{2}_{L^2 (\R^d)} + \|(I - P_{\lam})f\|^{2}_{L^2 (\R^d)}\\
&\leq C\exp(C\lam)\|P_{\lam}f\|^{2}_{L^{2}(\Omega)} + \|(I - P_{\lam})f\|^{2}_{L^2 (\R^d)}\\
&= C\exp(C\lam)\|f - (I - P_{\lam})f\|^{2}_{L^{2}(\Omega)} + \|(I - P_{\lam})f\|^{2}_{L^2 (\R^d)}\\
&\leq 2C\exp(C\lam)\|f\|^{2}_{L^2 (\Omega)} + 2C\exp(C\lam)\|(I-P_{\lam})f\|^{2}_{L^2 (\Omega)} + \|(I-P_{\lam})f\|^{2}_{L^2 (\R^d)}\\
&\leq 2C\exp(C\lam)\|f\|^{2}_{L^2 (\Omega)} + (2C\exp(C\lam) + 1)\|(I-P_{\lam})f\|^{2}_{L^2 (\R^d)} .
\end{align*}
Also, by Lemma 1 in \cite{green2020energy}, we have 
\begin{align*}
c \|(I-P_{\lam})f\|^{2}_{L^2 (\R^d)} \leq \|((-\Delta + 1)^{s/2} - \lam)f\|^{2}_{L^2 (\R^d)}
\end{align*}
for some $c > 0$ independent with $\lambda$.
Therefore, we conclude that
\begin{align*}
\|f\|^{2}_{L^2 (\R^d)} \leq C\exp(C\lam)\left[\|((-\Delta + 1)^{s/2} - \lam)f\|^{2}_{L^{2}(\R^d)} + \|f\|^{2}_{L^{2}(\Omega)}\right]. 
\end{align*}
\end{proof}
\begin{prop} \label{prop:resolvent2}
Let $s \geq 2$ and assume that $\Omega\subset\R^d$ is thick. Then there exist $C,c > 0$ such that for all $U = (u_{1},u_{2})\in \hshd$ and all $\lam\in\R$, 
\begin{align*}
c\exp(-C|\lam|)\|U\|^{2}_{\hsld} \leq \|(\Az - i\lam I)U\|^{2}_{\hsld} + \|u_{2}\|^{2}_{L^{2}(\Omega)}. 
\end{align*}
\end{prop}
\begin{proof}[Proof of Proposition \ref{prop:resolvent2}]
For $U = (u_{1},u_{2})\in\hshd$, we set 
\begin{align*}
\begin{pmatrix}
w_{1} \\
w_{2}
\end{pmatrix}
= 
\begin{pmatrix}
(-\Delta + 1)^{s/4} & -i\\
(-\Delta + 1)^{s/4} & i
\end{pmatrix}
\begin{pmatrix}
u_{1}\\
u_{2}
\end{pmatrix}. 
\end{align*}
By the parallelogram law, we obtain
\begin{align*}
\|w_{1}\|^{2}_{L^2 (\R^d)} + \|w_{2}\|^{2}_{L^2 (\R^d)} = 2\|U\|^{2}_{\hsld}. 
\end{align*}
Moreover, we have 
\begin{align*}
\|(\Az - i\lam I)U\|^{2}_{\hsl} &= \|(-\Delta + 1)^{s/2}(-i\lam u_{1} + u_{2})\|^{2}_{L^2} + \|-(-\Delta + 1)^{s/2}u_{1} - i\lam u_{2}\|^{2}_{L^2}\\
&= \|-\lam\frac{w_{1} + w_{2}}{2} + (-\Delta + 1)^{s/2}\frac{w_{1} - w_{2}}{2}\|^{2}_{L^2} \\
&\quad+ \|-(-\Delta + 1)^{s/2}\frac{w_{1} + w_{2}}{2} + \lam\frac{w_{1} - w_{2}}{2}\|^{2}_{L^2}\\
&=\|\lam w_{1} - (-\Delta + 1)^{s/2}w_{1}\|^{2}_{L^2} + \|\lam w_{2} + (-\Delta + 1)^{s/2}w_{2}\|^{2}_{L^2}. 
\end{align*}
For $\lam \geq 0$, applying Proposition \ref{prop:resovent} to $w_{1}$ with $s/2$, we have
\begin{align*}
&2c\exp(-C\lam)\|U\|^{2}_{\hsl} \\
&= c\exp(-C\lam)(\|w_{1}\|^{2}_{L^2} + \|w_{2}\|^{2}_{L^2})\\
&\leq \|((-\Delta + 1)^{s/4} -\lam)w_{1}\|^{2}_{L^2} + \|w_{1}\|^{2}_{L^2 (\Omega)} + c\exp(-C\lam)\|w_{2}\|^{2}_{L^2}\\
&\leq \|((-\Delta + 1)^{s/4} -\lam)w_{1}\|^{2}_{L^2} + 2\|w_{1} -w_{2}\|^{2}_{L^2 (\Omega)} + c\|w_{2}\|^{2}_{L^2}\\
&\leq \|((-\Delta + 1)^{s/4} -\lam)w_{1}\|^{2}_{L^2} + c\|((-\Delta + 1)^{s/4} +\lam)w_{2}\|^{2}_{L^2} + 8\|u_{2}\|^{2}_{L^2(\Omega)}\\
&\leq c\|(\mathcal{A}_{0} - i\lam I)U\|^{2}_{\hsl} + 8\|u_{2}\|^{2}_{L^2(\Omega)}. 
\end{align*}
For $\lam < 0$, we get the same inequality replacing the role of $w_{1}$ with $w_{2}$. 
\end{proof}
\subsection{Energy decay}
Finally we prove $\eqref{item:thickg} \implies \eqref{item:logenergy}$.
By the assumption \eqref{item:thickg}, $\Omega = \{\gamma \geq \varepsilon\}$ is thick for some $\varepsilon > 0$. 
Therefore, by Proposition \ref{prop:resolvent2}, we have
\begin{align*}
c\exp(-C|\lam|)\|U\|^{2}_{\hsl} &\leq \|(\Az - i\lam I)U\|^{2}_{\hsl} + \|u_{2}\|^{2}_{L^{2}(\Omega)}\\
&\leq 2\|(\Ag - i\lam I)U\|^{2}_{\hsl} + (2 + \varepsilon^{-2})\|\gamma u_{2}\|^{2}_{L^2 (\Omega)}. 
\end{align*}
Since $\Az$ is skew-adjoint, we obtain
\begin{align*}
\Re \lan(\Ag - i\lam I)U,U\ran = \Re \lan(\Az - i\lam I)U,U\ran - \lan\gamma u_{2},u_{2}\ran = -\|\sqrt{\gamma} u_{2}\|^{2}_{L^2}. 
\end{align*}
By the Cauchy--Schwarz inequality, we have 
\begin{align*}
D\|\gamma u_{2}\|^{2}_{L^2} \leq \|\gamma\|_{L^{\infty}}\|\sqrt{\gamma}u_{2}\|^{2}_{L^2}\leq \frac{D^2 \|\gamma\|^{2}_{L^{\infty}}\|(\Ag - i\lam)U\|^{2}_{\hsl}}{\delta} + \delta\|U\|^{2}_{\hsl}. 
\end{align*}
for any $D,\delta > 0$. 
Taking $D = 2 + \varepsilon^{-2}$ and $\delta = c\exp(-C|\lam|)/2$, we obtain
\begin{align*}
&c\exp(-C|\lam|)\|U\|^{2}_{\hsl} \\
&\leq 2\|(\Ag - i\lam I)U\|^{2}_{\hsl} + (2 + \varepsilon^{-2})\|\gamma u_{2}\|^{2}_{L^2 (\Omega)}\\
&\leq 2\|(\Ag - i\lam I)U\|^{2}_{\hsl} + \frac{(2 + \varepsilon^{-2})^2 \|\gamma\|^{2}_{L^{\infty}}}{c\exp(-C|\lam|)}\|(\Ag - i\lam I)U\|^{2}_{\hsl}\\
&\quad +\frac{1}{2}c\exp(-C|\lam|)\|U\|^{2}_{\hsl}. 
\end{align*}
By this inequality, we have
\begin{align*}
c\exp(-C|\lam|)\|U\|^{2}_{\hsl} &\leq \|(\Ag - i\lam I)U\|^{2}_{\hsl}, 
\end{align*}
here the constants $c, C$ may differ from the previous ones.
Applying Theorem \ref{thm:burq} with $k = 1$, we conclude that \eqref{item:logenergy} holds. 

        \section*{Acknowledgment}
	The authors would like to thank Professor Mitsuru Sugimoto for valuable discussions. 

\end{document}